\documentclass[11pt]{article}

\usepackage{amsxtra,amstext}
\usepackage{amssymb}
\usepackage{latexsym}

\newtheorem{theor1}{Theorem}
\newtheorem{theor3}{Proposition}
\newtheorem{theor4}{Lemma}
\newtheorem{theor5}{Corollary}

\addtocounter{footnote}{+2}

\title{On the additive theory of prime numbers II}
\vspace{0,5cm}

\author{Patrick {\sc Cegielski}\thanks{LACL, UMR--FRE 2673, Universit\'e
Paris 12, IUT Route Foresti\`ere Hurtault F-77300 Fontainebleau, \hfill\break
- Email: cegielski@univ-paris12.fr},
\ Denis {\sc
Richard}\thanks{LLAIC1 Universit\'e d'Auvergne, IUT Informatique,
B.P. 86, F-63172 Aubi\`ere Cedex \hfill \break
- Email: richard@iut.u-clermont1.fr}
\ \& Maxim {\sc
Vsemirnov}\thanks{Steklov Institute of Mathematics (POMI), 27 Fontanka St
Petersburg, 191011, Russia \hfill\break
- Email: vsemir@pdmi.ras.ru }}
\date{\today}
\vspace{3mm}

\begin{document}

\maketitle

\begin{abstract}
The undecidability of the additive theory of primes
(with identity) as well as the theory ${\rm Th} (\mathbb{N},+, n
\mapsto p_n)$, where $p_n$ denotes the $(n+1)$-th prime, are open
questions. As a possible approach, we extend the latter theory by
adding some extra function. In this direction we show the
undecidability of the existential part of the theory ${\rm Th}
(\mathbb{N}, +, n \mapsto p_n, n \mapsto r_n)$, where $r_n$ is the
remainder of $p_n$ divided by $n$ in the euclidian division.
\end{abstract}

\begin{quote}
\centerline{\bf R\'esum\'e}
{\small \quad 
L'ind\'ecidabilit\'e de la th\'eorie additive des nombres premiers ainsi que
de la th\'eorie ${\rm Th} (\mathbb{N},+, n \mapsto p_n)$, o\`u $p_n$ d\'esigne
le $(n+1)$-i\`eme premier, sont deux questions ouvertes. Nous \'etendons cette
derni\`ere th\'eorie en lui ajoutant une fonction suppl\'ementaire et nous
montrons l'ind\'ecidabilit\'e de la th\'eorie
${\rm Th} (\mathbb{N}, +, n \mapsto p_n, n \mapsto r_n)$, o\`u $r_n$ d\'esigne
le reste de $p_n$ de la division euclidienne de $p_n$ par $n$, et m\^eme de sa
seule partie existentielle.}
\end{quote}

\underline{\bf Introduction} - The additive theory of primes contains longtime open
classical conjectures of Number Theory, as famous {\sc Goldbach}'s binary
conjecture or {\sc twin primes} conjecture, and so on. Some authors provided
[BJW,BM,LM] conditional proofs (through {\sc Schinzel}'s Hypothesis [SS]) of
the undecidability of the additive theory of primes
${\rm Th} (\mathbb{N}, +, \mathbb{P})$, where $\mathbb{P}$ is the set of all
primes. Weakening the problem by strengthenning this theory, we introduced
[CRV] the theory ${\rm Th}(\mathbb{N},+, n\mapsto p_n)$, where $p_n$ is the
$(n+1)$-th prime, and posed the problem of its (un)decidability. As usual for
a language containing a function symbol, we suppose it contains identity. Note
that $\mathbb{P}$ is existentially definable within
$\langle \mathbb{N}, n \mapsto p_n \rangle$, hence
${\rm Th}(\mathbb{N},+, \mathbb{P})$ is a subtheory of
${\rm Th}(\mathbb{N},+, n \mapsto p_n)$. At the moment, the undecidability of
the latter theory is still an open question, and our approach in [CRV] was to
consider several approximations of the function $n \mapsto p_n$ as, for
instance, $n \left\lfloor\log n \right\rfloor$ and on this way we showed the
undecidability of theories ${\rm Th}(\mathbb{N}, +, nf(n)$) for a family of
functions $f$ including $\left\lfloor \log \right\rfloor$ mentioned above.
Another approach consists of extending the language $\{+,n \mapsto p_n \}$ to
$\{ +, n \mapsto p_n,$ $n \mapsto r_n \}$, where $r_n$ is the remainder of
$p_n$ divided by $n$. The main result of this paper is the following:

\begin{theor1} Multiplication is existentially $\rm \langle \mathbb{N},+,
n \mapsto p_n,$ $\rm n \mapsto r_n \rangle$-definable at first-order.
\end{theor1}

This leads to the following (without use of conjectures) result:

\begin{theor5} \ $\rm {\rm Th}_\exists( \mathbb{N},+, n \mapsto p_n,$
$\rm n \mapsto r_n)$ is undecidable.
\end{theor5}

\underline{Remark} Actually all positive integer constants are existentially
$\{+,\mathbb{P} \}$-definable in the following manner:
$$
\begin{array}{lll}
x=0 & \Leftrightarrow & x+x=x;\\
x=1 & \Leftrightarrow & \exists y (y=x+x \wedge y \in \mathbb{P});\\
x=2 & \Leftrightarrow & \exists y (y=1 \wedge x=y+y);\\
& ~\vdots &\\
x=n+1 & \Leftrightarrow & \exists y \exists z (y= n \wedge z=1 \wedge x=y+z).
\end{array}
$$
As we mentioned above, $\mathbb{P}$ is existentially definable within the
language $\{+, n \mapsto p_n \}$, hence all positive integer constants are
also existentially $\{+, n \mapsto p_n \}$-definable. Note, that
$n \left\lfloor \frac{p_n}{n} \right\rfloor = p_n - r_n$. We intend to define
(section 2, see Lemma 3) $\left\lfloor \frac{p_n}{n} \right\rfloor$ from + and
$n \left\lfloor \frac{p_n}{n} \right\rfloor$. Then the strategy will be to
define multiplication through the function $n \mapsto cn^2$ (where $c$ is a
fixed constant), which is to be proved $\{+, \left\lfloor \frac{p_n}{n}
\right\rfloor , n \left\lfloor \frac{p_n}{n} \right\rfloor\}$-definable.
Consequently, multiplication will be existentially
$\{+, n \mapsto p_n$, $n\mapsto r_n\}$-definable at first-order.

\

\underline{Remark.} In the previous paper [CRV] we consider continuous real
strictly increasing functions and their inverses. Since we work with integer
parts we have to introduce pseudo-inverses of discrete functions. For such a
discrete unbounded function $f$ from $\mathbb{N}$ into $\mathbb{N}$, we define
its pseudo-inverse $f^{-1}$ from $\mathbb{N}$ into $\mathbb{N}$ by
$f^{-1}(n)= \mu m[f(m+1) >n]$, where $\mu$ means ``the smallest $\ldots$ such
that''. Due to the unboundness of $f$ such an $f^{-1}$ is correctly defined.

\

{\large \bf 1) Some preliminary results in Number Theory}

\

Contrarily to what happens with $\log$, the behavior of
$\left\lfloor\frac{p_n}{n} \right\rfloor$ is {\it a priori} irregular but we
shall prove it is not too much chaotic. Namely, we prove:

\begin{theor3} \label{P11} \
Let us denote the mapping $n \mapsto \left\lfloor \frac{p_n}{n} \right\rfloor$
by $f$.

1) For $m>n$, we have $f(m) - f(n) \geq -1$;

2) For $n \geq 11$, we have $f^{-1} (n+1) - f^{-1}(n) >n$.
\end{theor3}

{\bf Proof} \hspace{0.2cm}
1) We use the following estimates for $p_n$ ([RP], p. 249):

\medskip

$p_m \geq m\log m + m \log \log m - 1.0072629 m$ \ for $m \geq 2$;

\smallskip

$p_m \leq m \log m + m \log \log m - 0.9385 m$ \ for $m \geq 7022$. \hfill (1)

\medskip

For $m > n \geq 7022$, we have $f(m) - f(n) = \left\lfloor \frac{p_m}{m}
\right\rfloor - \left\lfloor \frac{p_n}{n} \right\rfloor$

$ \geq \frac{p_m}{m} - \frac{p_n}{n} - 1
\geq log(\frac{m}{n}) - log(\frac{log m}{log n}) - 0.9385 + 1.0072629 - 1$.

Hence $f(m) - f(n) \geq -1$ because the sum of the two first terms is
positive as is the sum of terms three and four.

If $n < 7022$, one may check the desired inequality by a direct computation.

\medskip

2) Let $m$ be $f^{-1}(n)$. By the very definition of $f^{-1}$, the equality
$m=f^{-1}(n)$ is equivalent to the conjunction of the two following conditions:

$$
\hspace{5,5cm}\left \{
\begin{array}{lll}
\displaystyle\left\lfloor \frac{p_{m+1}}{m+1} \right\rfloor \geq n+1;\\
\displaystyle\forall k \leq m \ \left\lfloor \frac{p_k}{k}
\right\rfloor \leq n.\hspace{5cm}(2)
\end{array}
\right.
$$

\newpage

For $k \leq 7022$, the maximum of $\frac{p_k}{k}$ is attained for $k=7012$ and
equal to $\frac{p_{7012}}{7012} < 10.102824 < 11$. Consequently, we see that
$m = f^{-1} (n) \geq f^{-1}(11) \geq 7022$ and this is the reason why in the
hypothesis of Proposition 1, item 2) we assume $n \geq 11$.

To prove the inequality, it is sufficient to prove that for $k = m+n$ we
have $\left\lfloor \frac{p_k}{k} \right\rfloor \leq n+1$, or in other words,
$$\frac{p_k}{k} < n+2. \eqno \rm{(3)}$$

Note that for $m \geq 7022$, we have by (2):

$$n+1 \leq \left\lfloor \frac{p_{m+1}}{m+1} \right\rfloor +1 \leq
\frac{p_{m+1}}{m+1}+1 \leq \log (m+1) +\log\log(m+1) - 0.07 < m.$$
Consequently it is sufficient -- and actually more convenient -- to prove a
somehow stronger result, namely the same inequality (3) but for $m \geq 7022$
and $m+1 \leq k \leq 2m.$

From the second estimate of (1) we have, since
$k \geq m \geq 7022$, the following inequalities:
$$
\begin{array}{lll}
\frac{p_k}{k} &<& \log k + \log\log k - 0.9385\\
&\leq &\log 2m + \log\log 2m - 0.9385\\
&=& \log m +\log\log m + \log 2 + \log (1+\frac{\log 2}{\log m}) - 0.9385;
\end{array}
$$
using the first estimate of (1) and
$\frac{\log 2}{\log m}\leq \frac{\log 2}{\log 7022}$, we have:
$$\log m +\log \log m-1.0072629 \leq \frac{p_m}{m};$$
consequently:
$$\frac{p_k}{k} \leq \frac{p_m}{m} + 0.07 + \log 2 +
\log (1+ \frac{\log 2}{\log 7022}) \leq \frac{p_m}{m} +1$$
by an easy computation and finally, due to (2), we obtain
$\frac{p_k}{k} < n+2.$ \hfill $\Box$

\

Item 1) of previous proposition emphasizes on the fact that
$f : n \mapsto \left\lfloor \frac{p_n}{n} \right\rfloor$ is ``almost''
increasing and item 2) shows that the difference $f^{-1}(n+1) - f^{-1}(n)$ is
big enough with respect to $n$. This suggests to introduce a new class of
functions, containing $f$, for which we prove that the existential part of
the theory ${\rm Th} (\mathbb{N} , +, n \mapsto nf(n))$ is undecidable.

\

{\large \bf 2) The class $C(k,d,n_0)$ and some its properties}

\medskip

{\bf 2.1) \ The class $C(k,d,n_0)$}

\medskip

Let $k \geq 0$ be a fixed nonnegative integer. We shall say $f$ is
$k$-{\it almost increasing} if and only if
$$\forall y \geq x [f(y)- f(x) \geq - k]. \eqno {\rm(4)}$$

\medskip

In this sense $0$-almost increasing means increasing (not necessarily
strictly) and $n \mapsto \left\lfloor\frac{p_n}{n} \right\rfloor$ is
$1$-almost increasing (due to Proposition 1).

\medskip

Still looking at $n \mapsto \left\lfloor\frac{p_n}{n} \right\rfloor$, we
intend to consider functions whose pseudo-inverse is defined and
asymptotically increases quickly enough with respect to its argument. Let us
say that $f^{-1}$ has at least $(1/d)$-linear difference, if
$$\exists n_0 \in \mathbb{N} \forall n \geq n_0 [f^{-1} (n+1)-f^{-1}(n) >
\frac{n}{d}]. \eqno {\rm (5)}$$
In fact, for $\left\lfloor \frac{p_n}{n} \right\rfloor$, the constant $d$
is 1 and $n_0 = 11$, but results and proofs hold for an arbitrary (fixed) $d$.
\medskip

Now let us definite the class $C(k,d,n_0)$ as the set of functions from
$\mathbb{N}$ into $\mathbb{N}$ satisfying conditions (4) of being
$k$-{\it almost increasing} and (5) of having its pseudo-inverse with an
{\it at least} $(1/d)$-{\it linear difference}.

\smallskip

In order to prove {\sc fundamental lemma} of section 3, whose Theorem 1 is a
corollary, we show some properties of the class $C(k,d,n_0)$. Firstly, in
section 2.2 we present in three lemmas these properties and comment them.
Afterwards, in section 2.3, we prove them.

\

{\bf 2.2) \ Properties of $C(k,d,n_0)$}

\begin{theor4} \label{L21} For any function $f \in C(k,d,n_0)$ the following
items hold:

\smallskip

\quad (i) For any $n \geq n_0$, we have $f^{-1} (n+d) - f^{-1} (n) >n$;

\smallskip

\quad (ii) For any $n \geq n_0+1$, the set
$\{ x \in \mathbb{N} \; | \; f(x) =n \}$ is nonempty;

\smallskip

\quad (iii) For any $n \geq n_0+1$, the equality $f(x)=n$ implies
$$x > \frac{1}{2d} [(n-1)(n-2) -n_0(n_0-1)].$$
\end{theor4}

\medskip

\begin{theor4} \label{L22} If $f \in C(k,d,n_0)$ and $f(x) = n \geq n_0$, then
for any $c$ such that $1 \leq c \leq n$, we have:
$$-k \leq f(x+c)-f(x) \leq k+d. \eqno {\rm (6)}$$
\end{theor4}

\medskip

\begin{theor4}\label{L23} For any $f \in C(k,d,n_0)$, let
$x_0 = f^{-1} (2+4d+n_0^2 +k).$

Consider
$\tilde{f} : [x_0 +1, + \infty [ \cap \mathbb{N} \longrightarrow \mathbb{N}$
with $\tilde{f} (x) = f(x)$. Then $\tilde{f}$ is existentially definable at
first-order within $\langle \mathbb{N},+,1,x \mapsto xf(x) \rangle$.
\end{theor4}

\underline{Remarks} 1) Item (i) of Lemma \ref{L21} provides a linear lower
bound of values of $f^{-1}$ when difference of arguments is the parameter $d$
of the considered class.

Item (ii) of the same lemma insure that $f$ is asymptotically onto, and item
(iii) gives a quadradic lower bound for solutions of the equation $f(x) =n$,
that we need in section 3.

\medskip

2) Actually, as the reader will see within the proof, Lemma \ref{L21} does not
use condition (4) of being $k$-almost increasing.

\medskip

3) Lemma \ref{L22} provides asymptotical bounds for the difference of two
values of $f$ with arguments taken in a short interval with respect to the
values of these arguments. Refering to the previous Lemma \ref{L21} we see
that $n$ is at most $\sqrt{2dx+ n_0^2}+2$.

\medskip

4) Lemma \ref{L23} generalizes the situation of the main result of the previous
paper [CRV] of the same authors when $\left\lfloor \log n \right\rfloor$ was
``extracted'', {\it i.e.} defined, from + and
$n \left\lfloor \log n \right\rfloor$.

\

{\bf 2.3) \ Proofs of the three Lemmas}

\

{\bf Proof of Lemma \ref{L21}} \ (i) By condition (5):
$$
\begin{array}{lll}
f^{-1}(n+d) - f^{-1}(n) &=& [f^{-1}(n+d) - f^{-1}(n+d-1)]\\
                        && + [f^{-1}(n+d-1) - f^{-1}(n+d-2)]\\
                        && + \ldots \\
                        && + [f^{-1}(n+1) - f^{-1}(n)]\\
                        &>& \frac{n+d-1}{d} + \frac{n+d-2}{d} + \ldots
                            + \frac{n}{d} \\
                        &>& n.
\end{array}
$$

(ii) If there was no $x$ such that $f(x) =n$, we would have $f^{-1}
(n)=f^{-1}(n-1)$. But $f^{-1}(n) > f^{-1} (n-1)$ according to condition (5).

\medskip

(iii) By definition of $f^{-1)}$, we have: $x > f^{-1}(n-1)$.

As in (i), we have:
$$
\begin{array}{lll}
f^{-1}(n-1) - f^{-1}(n_0) &=& [f^{-1}(n-1) - f^{-1}(n-2)]\\
                        && + \ldots \\
                        && + [f^{-1}(n_0+1) - f^{-1}(n_0)]\\
                        &>& \frac{n-2}{d} + \frac{n_0}{d} + \ldots
                            + \frac{n}{d} \\
                        &=& \frac{(n-2)(n-1) - n_0(n_0+1)}{2d}.
\end{array}
$$
and the result. \hfill$\Box$

\medskip

{\bf Proof of Lemma \ref{L22}} \ The left-hand side of the inequality is an immediate
consequence of the very definition of a $k$-almost increasing function. For
proving the right-hand side, note that, using $k$-almost increasing property
of $f$ together with $f(x) =n$, we obtain:
$$\displaystyle\max_{y \leq x} f(y) \leq f(x) + k = n+k,$$
so that $f^{-1}(n+k) \geq x$, by the definition of $f^{-1}$. By previous Lemma
\ref{L21}, item (i) and the latter inequality, we have:
$$f^{-1}(n+k+d)>f^{-1}(n+k)+n+k \geq x + n + k \geq x+n \geq x+c$$
since $1\leq c\leq n$. Using again the definition of $f^{-1}$, we see that
$f(x+c)\leq n+k+d=f(x)+k+d$ and we are done. \hfill $\Box$

\bigskip

{\bf Proof of Lemma \ref{L23}} \ To define $\tilde{f}$ within the structure
$\langle \mathbb{N}, +,x \mapsto xf(x) \rangle$ we shall make use of the
inequality:
$$0 \leq f(x) < x$$
together with the remainder of $f(x)$ modulo $x+1$, which we must define in the
considered structure.

\medskip

{\bf Fact 1}.- $f(x) < x$.

\medskip

By the definition of $f^{-1}$, we have $f (x_0+1)>k+2+4d+n_0^2$ and by the
$k$-almost increasing property we deduce, for $x \geq x_0 +1$,
$$n=f(x) \geq f(x_0+1)-k > 2+4d+n_0^2. \eqno {\rm (7)}$$
Hence $\frac{n-2}{2d} >2$.

From (7), we obtain $n>n_0+1$ so that by Lemma \ref{L21}, item (iii), we have:

$$x>\frac{1}{2d}[(n-1)(n-2)-n_0(n_0-1)],$$

\medskip 

hence:

$$x > 2(n-1) - \frac{n_0(n_0-1)}{2d}
    > 2(n-1) - n_0^2 = n+ (n- 2 - n_0^2) > n = f(x). \; \Box\Box$$

\medskip

{\bf Fact2}.- {\em We have}:
$$f(x) \equiv (x+1)f(x+1) - xf(x) \hspace{0.2cm} (\bmod \, x+1).\eqno
{\rm (8)}$$

It is sufficient to note that
$(x+1)f(x+1)-xf(x)=f(x)+(x+1)[f(x+1)-f(x)]. \; \Box\Box$

\medskip

We are still unable to define general congruences, fortunately here the
difference \break $|f(x+1)-f(x)|$ is bounded, namely,
$$|f(x+1) -f(x)| \leq k+d, \eqno {\rm (9)}$$
due to Lemme \ref{L22}, with $c=1$. This suggests to introduce the notion of a
restricted congruence, namely, for $a,b,m$ in $\mathbb{N}$ and some fixed
integer $c$, we define $a \equiv_c b (\bmod \, m)$ by:
$$\bigvee^c_{h=0} \{[a =b +{\underbrace{m+\cdots +m}_{\displaystyle{\rm h
\ times}}} ] \vee [ b=a + {\underbrace{m+\cdots + m}_{\displaystyle{\rm
h \ times}}}]\}.$$

Obviously, the first-order latter formula is expressible within the structure
$\langle \mathbb{N},+ \rangle$, since $c$ is fixed. The congruence (8) and
inequality (9) provide together the following restricted congruence:
$$f(x) \equiv_{k+d} (x+1)f(x+1)-xf(x)(\bmod \, x+1),$$

which is a definition of
$f(x)$ within $\langle \mathbb{N}, +, 1, x \mapsto xf(x) \rangle$
since $1\leq f(x) < x$. Finally, we provide explicitely an existential
first-order definition of $f$, namely:
$$[x > x_0 \wedge y = f(x)] \leftrightarrow$$

$$ \{x>x_0 \wedge y \leq x \wedge
\bigvee^{k+d}_{h=0} [(y + xf(x)=
(x + 1)f(x + 1)+{\underbrace{(x + 1) + \cdots + (x + 1)}_{\displaystyle{\rm h
\ times}}})$$

$$\vee((x+1)f(x+1)=y+xf(x)+{\underbrace{(x+1)+ \cdots +
(x+1)}_{\displaystyle{\rm h \ times}}})] \}.$$

\

{\large \bf 3) \ Fundamental Lemma and the proof of the Main Theorem}

\medskip

In order to prove the undecidability of
${\rm Th}(\mathbb{N},n \mapsto p_n, n \mapsto r_n)$, we prove a more general
result, namely:

\begin{theor4} [Fundamental Lemma] \label{L31} \ For any $f \in C(k,d,n_0)$
{\rm [see \S 2]}, multiplication is existentially
$\{ +, \; 1, \; x \mapsto xf(x) \}$-definable at first-order.
\end{theor4}

As shown by Y. {\sc Matiyasevich}, the existential true theory of numbers is
exactly the set of arithmetical relations, which are definable by diophantine
equations. Therefore the negative solution of the 10-th Hilbert's problem [MY]
implies the following corollary.

\begin{theor5} \ The existential theory
${\rm Th}_\exists (\mathbb{N}, \; +, \;1, \; x \mapsto xf(x))$ is undecidable.
\end{theor5}

{\bf Proof of Lemma \ref{L31}} \ It suffices to show that for some constants $c$ and
$n_1$ the function $n \mapsto cn^2$ from $[n_1, + \infty [ \cap \mathbb{N}$
into $\mathbb{N}$ is $\{+, \; 1, \; x \mapsto xf(x) \}$-definable. More
precisely, we
shall take $c=5d$ and $n_1=2+5d+n_0^2$. Consider $n \geq n_1$. Since
$n_1 > n_0+1$, we can apply Lemma \ref{L21}, item (ii), proving there exists
$x$ such that $f(x) =5dn$. Let $x_0$ be the same as in Lemma \ref{L23}, namely
$x_0 =f^{-1} (2+4d+n_0^2 +k)$. Let us show $x > x_0$. Otherwise $x \leq x_0$,
so that by the $k$-almost increasing property $f(x) \leq f(x_0)-k$, implying,
by the definitions of $f^{-1}$ and $x_0$,
$$f(x) \leq 2+4d+n_0^2+k-k<n_1<5dn_1 \leq 5dn=f(x),$$
which is impossible.

Note that $5dn$ is $\{+ \}$-definable as the sum of $5d$ terms equal to $n$
($d$ is a fixed constant). Now thanks to Lemma \ref{L23}, an $x$ such that
$f(x)= 5dn$ is $\{+,1,x \mapsto xf(x) \}$-definable.

On the other hand:
$$(x+n)f(x+n)-xf(x)=(x+n)[f(x+n)-f(x)]+nf(x) =(x+n)[f(x+n)-f(x)]+5dn^2.$$
By Lemma \ref{L22} applied to $c=n$, we have $|f(x+n)-f(x)| \leq k+d$, so that:
$$5dn^2 \equiv_{k+d} (x+n)f(x+n)-xf(x)(\bmod \, x+n). \eqno {\rm (10)}$$

According to Lemma \ref{L21} and item (iii) since $f(x) =5dn$ and
$5dn >n_1 > n_0+1$ the inequalities $n \geq n_1 > n_0^2$ and:

\begin{eqnarray}
\hspace{3cm}x+n &> & \frac{(5dn-1)(5dn-2)}{2d} - \frac{n_0(n_0-1)}{2d}
+n\nonumber\\
\hspace{3cm} &>& \frac{25d^2n^2 -15nd}{2d} > 5dn^2 \nonumber
\hspace{5.5cm} {\rm (11)}
\end{eqnarray}
hold.

Using (10) and (11), a similar argument as in Lemma \ref{L23} shows that the
function $n \mapsto 5dn^2=cn^2$ with domain $[n_1, + \infty [\cap \mathbb{N}$
is existentially $\{ +,1,x \mapsto xf(x)\}$-definable. By a routine argument,
multiplication is clearly existentially $\{+,1,x \mapsto xf(x)\}$-definable.
\hfill $\Box$

\

{\bf Proof of the Main-Theorem} \ We remind the reader that 1 was existentially
$\{+, \; \mathbb{P} \}$ and $\{+, \; n \mapsto p_n \}$-defined in the
introduction.

We also noted that $n \left\lfloor \frac{p_n}{n} \right\rfloor =p_n -r_n$ and
$n \mapsto n \left\lfloor \frac{p_n}{n} \right\rfloor$ belongs to $C(1,1,11)$,
the latter due to Proposition \ref{P11}, \S 1. Then Fundamental Lemma can be
applied and multiplication is existentially
$\{+,n \mapsto p_n, n \mapsto r_n \}$-definable. \hfill $\Box$

\

{\bf Conclusion}: Our main result is absolute in the sense that does not depend on
any conjecture. In order to shed more light on the considered theories
${\rm Th}_\exists (\mathbb{N},+,\mathbb{P})$ and
${\rm Th}_\exists (\mathbb{N}, n \mapsto p_n, n \mapsto r_n)$, we recall a
conditional result of A. {\sc Woods}. Let us recall that {\sc Dickson}'s
{\sc conjecture} [DL] claims that if
$a_1, a_2, \ldots a_n, b_1,b_2, \ldots b_n$ are integers with all $a_i > 0$ and
$$\forall y \neq 1 \exists x [y \not\kern 3pt\mid \prod_{1\leq i\leq n}
(a_i x + b_i)]$$
then there exist infinitely many $x$ such that $a_i x +b_i$ are primes for all
$i$. Let us call $DC$ this conjecture, then A. {\sc Woods} proved [WA]:

\smallskip

{\em If $DC$ is true then the existential theory
${\rm Th}_\exists (\mathbb{N},+,\mathbb{P})$ is decidable.}

\smallskip

Now, the question is to know whether there is a gap between
${\rm Th}_\exists (\mathbb{N},+,n \mapsto p_n,n \mapsto r_n)$ or whether they
are exactly the same. In the case of equality between these two theories, $DC$
is false (and hence {\sc Schinzel}'s {\sc hypothesis} on primes, whose $DC$ is
the linear case, is also false).

\

{\bf Open problem:} \ {\em Is ${\rm Th}_\exists (\mathbb{N},+,\mathbb{P})$ equal to
${\rm Th}_\exists (\mathbb{N},+,n \mapsto p_n, n \mapsto r_n) ?$}

\

\centerline{\bf References}

\medskip

[BJW] \ P.T. {\sc Bateman}, C.G. {\sc Jockusch} and A.R. {\sc Woods},
{\it Decidability and Undecidability of theories with a predicate for the prime},
{\bf Journal of Symbolic Logic}, vol. 58, 1993, pp.672-687.

\medskip

[BM] \ Maurice {\sc Boffa}, {\it More on an undecidability result of {\sc Bateman,
Jockusch} and {\sc Woods}}, {\bf Journal of Symbolic Logic}, vol. 63, 1998, p.50.

\medskip

[CRV] \ Patrick {\sc Cegielski}, Denis {\sc Richard} \& Maxim {\sc Vsemirnov},
{\it On the additive Theory of Prime Numbers I}, {\bf Proceedings of CSIT'2003} (Computer
Science and Information Technologies), September 22-26, 2003, Yerevan,
Armenia, 459 p., pp. 80--85.

\medskip

[DL] \ L.E. {\sc Dickson}, {\it A new extension of {\sc Dirichlet}'s theorem on
prime numbers}, {\bf Messenger of Mathematics}, vol. 33 (1903--04), pp. 155--161.

\medskip

[LM] \ T. {\sc Lavendhomme} \& A. {\sc Maes}, {\it Note on the undecidability of}
$\langle \omega,+,P_{m,r} \rangle$, Definability in arithmetics and
computability, 61-68. {\bf Cahier du Centre de logique}, Belgium, 11 (2000).

\medskip

[MY] \ Yuri {\sc Matiyasevich}, {\bf  Hilbert's tenth Problem}, The MIT
Press, Foundations of computing, 1993, XXII+262p.

\medskip

[RP] \ Paul {\sc Ribenboim}, {\bf The new book of Prime records}, Springer, 1996,
XIV+541p.

\medskip

[SS] \ A. {\sc Schinzel} \& W. {\sc Sierpie\'{n}sky}, {\it Sur certaines
hypothèses concernant les nombres premiers}, {\bf Acta Arithmetica}, vol. 4, 1958,
185--208 and 5, 1959, 259.

\medskip

[WA] \ Alan {\sc Woods}, {\bf Some problems in logic and number theory, and their connection}, Ph.D.
thesis, University of Manchester, Manchester, 1981.

\end{document}